\documentstyle[11pt]{article}
                      
\title{{\Large SPECHT MODULES FOR FINITE REFLECTION GROUPS}} 

\author{$by$ S.~HALICIO{\u G}LU} 
\date{ }
\begin{document}
\maketitle

\noindent
{\bf 1. Introduction}

Over fields of characteristic zero, there are well known constructions
of the irreducible representations, due to A Young,  and of irreducible 
modules, called Specht modules, due to W Specht, for the symmetric groups 
$S_{n}$ 
which are based on elegant combinatorial 
concepts connected with Young tableaux etc.(see, e.g.{\bf [13]}).  
James {\bf [12]} extended these ideas to construct irreducible 
representations and modules over arbitrary field. Al-Aamily,   
Morris and Peel {\bf [1]} showed how this 
construction could be extended cover the  Weyl groups of type 
$B_{n}$. In {\bf [14]} Morris described a possible extension of 
James' work for Weyl groups in general. Later, the present author and
 Morris {\bf [8]}  give an alternative generalisation of James' work which 
is an extended improvement and extension of the original 
approach suggested by Morris. We now give a possible extension of James' work
 for finite reflection groups in general.

\noindent
{\bf 2. Some General Results On Finite Reflection Groups}

In this section we establish the notation and state some results on 
finite reflection groups which are required later . Standard references 
for this material
 are N Bourbaki {\bf [3]},  R W Carter {\bf [4]},   
 J E Humphreys {\bf [10] [11]}, Grove and
  Benson {\bf [7]}.

Let  $V$ be $l$-dimensional Euclidean space over the real field {\bf R}
 equipped with a positive definite inner product ( , ). For $\alpha \in V$,
 $\alpha \neq 0$, let $\tau_{\alpha}$ be the $reflection$ in the   
  hyperplane orthogonal to $\alpha$, that is, $\tau_{\alpha}$ is the
   linear transformation on $V$ defined by 
   \[\tau_{\alpha}(v)= v - 2\frac{(\alpha,v)}{(\alpha,\alpha)}\alpha\]
\noindent
for all $v \in V$. Let $\Phi$ be a root system in  $V$ and 
$\pi$  be a simple system in $\Phi$ with corresponding positive system 
$\Phi^{+}$ and negative system $\Phi^{-}$. Then the finite reflection
 group
\[ {\cal W}={\cal W}(\Phi)=<\tau_{\alpha}\mid \tau_{\alpha}^{2} = e, 
~~(\tau_{\alpha}\tau_{\beta})^{m_{\alpha\beta}}= e,~~ \alpha,\beta \in \pi~and~
\alpha \neq \beta  >\]
\noindent
where $e$ is the identity element of ${\cal W}$ and $m_{\alpha \beta}$ is the
 order of $\tau_{\alpha}\tau_{\beta}$. Let $l(w)$ denote the $length$ of $w$ 
and the $sign$ of $w$, 
$s(w)$, is defined by 
$s(w)=(~-1~)^{l(w)}~,~w \in {\cal W}$.

We note the following facts which are required later.

\noindent
 2.1. There are $\mid {\cal W} \mid$ simple systems (positive systems) 
in $\Phi$ 
given by $w \pi $ $( w \Phi^{+})$, $w \in {\cal W}$. The group ${\cal W}$ 
acts transitively on the set of simple systems.

\noindent
 2.2.  To each root system $\Phi$, there corresponds a graph $\Gamma$ 
called the $Coxeter~graph$ (or $Dynkin$ $diagram$) of 
${\cal W}$, whose nodes are $1:1$ correspondence with the elements of
 $\pi$. A finite reflection group is $irreducible$ if its Coxeter graph is 
 connected. 
Finite irreducible reflection groups have been classified and correspond to 
root systems of type $A_{l}(l\geq 1)$, $ B_{l} (l\geq 2) $, 
$C_{l} (l\geq 3)$, $ D_{l} (l \geq 4)$, $ ~E_{6}~, ~E_{7}~,~ E_{8} ~$,
$~F_{4}~,~G_{2} $, $H_{3}$, $H_{4}$, $I_{2}(~p~)~ (p = 5~ or~ p \geq 7)$. 
For example $ {\cal W}~(~A_{l}~)~\cong~S_{l+1} $, 
the symmetric group on the set $\{ 1 , 2 , ... , l + 1  \}$ . 
As our aim in this paper is to generalise ideas from the symmetric groups, 
the root system and simple system are given in this case. 
If $\{~\epsilon_{1}~,~\epsilon_{2}~,~...~,~\epsilon_{l+1}~\}$ is the standard basis for $R^{l+1}$ , then 
\[\pi ~=~\{~\alpha_{1}~=~\epsilon_{1}-\epsilon_{2}~,~\alpha_{2}~=~\epsilon_{2}-\epsilon_{3}~,~...~,~\alpha_{l}~=~\epsilon_{l}- \epsilon_{l+1}~\}\]  
\[\Phi~=~\{~\epsilon_{i}~-~\epsilon_{j}~\mid~1 \leq i ~,~j\leq l+1\}\]
\[\Phi^{+}~=~\{~\epsilon_{i}~-~\epsilon_{j}~\mid~1 \leq i ~<~j\leq l+1\}\]

\noindent
 2.3.  A $subsystem$ $\Psi$ of $\Phi$ is a subset of $\Phi$ which 
is itself a root system in the space which it spans. 
A subsystem $\Psi$ is said to be $additively~closed$ if $\alpha~,~\beta \in \Psi~,~ \alpha~+~\beta \in \Phi~$, then $ \alpha~+~\beta \in \Psi~$ . 
From now on we assume that $\Psi$ is additively closed 
subsystem of $\Phi$. A $finite~reflection~subgroup$ ${\cal W}(\Psi)$ of ${\cal W}$ corresponding to a 
subsystem $\Psi$ is the subgroup of ${\cal W}$ generated by the 
$\tau_{\alpha} ~,~\alpha \in \Psi $. If $\Psi$ and $\Upsilon$ are subsystems 
of $\Phi$ which span subspaces 
$U$ and $W$ of $V$ respectively , then $\Psi$ and $\Upsilon$ are $isomorphic$ if there exists a vector space isomorphism $\theta~:~U~\rightarrow W$ such that $\theta~(\Psi)=\Upsilon$ and 

\[~~\frac{~(~\theta~(~\alpha~)~,~\theta~(~\beta~)~)~}{~(~\theta~(~\beta~)~,~\theta~(~\beta~)~)~}=~\frac{~(\alpha~,~\beta~)}{(~\beta~,~\beta~)}~~~~~~~~~~~~~for ~all~ \alpha~,~\beta \in \Psi~.\]

\noindent
It follows that

\[{\cal W}~(\Psi)\cong {\cal W}~(\Upsilon)=\theta~{\cal W}~(\Psi)~\theta^{-1}~.\]

\noindent
The subsytems $\Psi$ and $\Upsilon$ are $conjugate~under~{\cal W}$ if $\Upsilon=w\Psi$ for some $w \in {\cal W}$ ; in this case ${\cal W}~(w~\Psi)=w~{\cal W}~(\Psi)~w^{-1}$ since $\tau_{w~(\alpha)}=w ~\tau_{\alpha} ~w^{-1}$
for $\alpha \in \Psi$ . Note that isomorphic subsystems are not necessarily conjugate.

\noindent
 2.4.  The graphs which are Dynkin diagrams of subsystems 
of $\Phi$ may be obtained up to conjugacy by a standard 
algorithm due independently to E B Dynkin, A Borel and J de Siebenthal (see
 e.g. {\bf [4]}).

\noindent
2.5.   If $w \in {\cal W}$ and 
$U$ is the subspace of $V$ composed of all vectors fixed by $w$, 
then $w$ is a product of reflections corresponding to roots in 
the orthogonal complement $U^{\perp}$ of $U$. {\bf [4]}

\noindent
2.6.  The simple system $J$ of $\Psi$ can always be chosen  such that $J\subset 
\Phi^{+}$. {\bf [15]}

\noindent
2.7.  The set $D_{\Psi}=\{ w \in {\cal W} \mid w (j) \in \Phi^{+}~ 
for~all~j \in J~\}$ is a 
$distinguished$ $ set$  $of$ $coset$ $representatives$ of ${\cal W}(\Psi)$ in ${\cal W}$ , 
that is, each element $w \in {\cal W}$ has unique expression of 
the form $d_{\Psi}w_{\Psi}$, where $d_{\Psi} \in D_{\Psi}$ and 
$w_{\Psi} \in {\cal W}(\Psi)$ and furthermore $d_{\Psi}$ is 
the element of minimal length in the coset $d_{\Psi}{\cal W}(\Psi)$. {\bf [15]}

\noindent
{\bf 3. Specht Modules for Finite Reflection Groups}

Let $\Phi$ be a root system with simple system $\pi$ and 
Coxeter graph $\Gamma$
and let $\Psi$ be a subsystem of $\Phi$ with simple system $J \subset \Phi^{+}$
 and  Coxeter graph $\Delta$.
If  $\Psi= \displaystyle \bigcup_{i=1}^{r}\Psi_{i} $ , 
where $\Psi_{i}$ are the indecomposable components of $\Psi$ , 
then let  $J_{i}$ be a simple system in $\Psi_{i}$ ($i=1,2,...,r$) 
and $J = \displaystyle \bigcup_{i=1}^{r}J_{i}$. Let $\Psi^{\perp}$ be the  largest
subsystem in $\Phi$ orthogonal to $\Psi$ and let $J^{\perp} \subset 
\Phi ^{+}$ the simple system of $\Psi^{\perp}$ .

Let  $\Psi^{'}$ be a subsystem of $\Phi$  which is contained in $\Phi \setminus \Psi$ , with simple system $J^{'} \subset \Phi^{+}$ and  
Coxeter graph $\Delta^{'}$. If  $\Psi^{'}=\displaystyle \bigcup_{i=1}^{s}\Psi_{i}^{'} $ , 
where $\Psi_{i}^{'}$ are the indecomposable components of $\Psi^{'}$  
then let  $J_{i}^{'}$ be a simple system in $\Psi_{i}^{'}$ ($i=1,2,...,s$) 
and $J = \displaystyle \bigcup_{i=1}^{s}J_{i}^{'}$. Let $\Psi^{'^{\perp}}$ be the  largest
subsystem in $\Phi$ orthogonal to $\Psi^{'}$ and let $J^{'^{\perp}} 
\subset \Phi ^{+}$ the simple system of $\Psi^{'^{\perp}}$ .

Let $\bar{J}$ stand for the ordered set $\{J_{1},J_{2},...,
J_{r};J_{1}^{'},J_{2}^{'},...,J_{s}^{'}\}$  , 
where in addition the elements in each 
$J_{i}$ and $J_{i}^{'}$ are ordered. Let 
\[{\cal T}_{J,J^{'}}=\{ w\bar{J} \mid w\in {\cal W} \} \]

Now , we  consider under what conditions the elements in the set 
${\cal T}_{J,J^{'}}$ are distinct. Such a condition is now obtained in 
the following lemma.

\noindent
{\sc 3.1. Lemma.} {\it  $\mid {\cal T}_{J,J^{'}}\mid=\mid {\cal W} \mid$ if and only if ${\cal W}(J^{\perp}) \cap {\cal W}(J^{'^{\perp}})=<~e~>$ .}

\noindent
{\it Proof.} See Lemma 3.1 {\bf [ 8 ]} $~~~~~~\Box$ 

Now we can give our principal definition .

\noindent
{\sc 3.2. Definition.} Let $\Psi$ and $\Psi^{'}$ be subsystems of $\Phi$ with 
simple systems $J$ and $J^{'}$ respectively such that 
$\Psi^{'} \subseteq \Phi \setminus \Psi$ and $J \subset \Phi^{+}$, 
$J^{'} \subset \Phi^{+}$. The pair $\{J,J^{'}\}$ is called a 
$useful ~system $ in $\Phi$ if ${\cal W}(J) \cap {\cal W}(J^{'})=<e>$ and 
${\cal W}(J^{\perp}) \cap {\cal W}(J^{'^{\perp}})=
<e>$.

\noindent
{\sc Remark 1.} If $\{J,J^{'}\}$ is a useful system in $\Phi$ , 
then $\{wJ,wJ^{'}\}$ is also a useful system in $\Phi$, 
for $w \in {\cal W}$. Thus , from now on ${\cal T}_{J,J^{'}}$ will be 
denoted by ${\cal T}_{\Delta}$. 

\noindent
{\sc Remark 2.} If $\{J,J^{'}\}$ is a useful system in $\Phi$ 
then $\Psi \cap \Psi^{'} = \emptyset$ and $\Psi^{\perp} \cap 
\Psi^{'^{\perp}} = \emptyset$. 
However the converse is not true in general. 

\noindent
{\sc 3.3. Definition.} Let $\{J,J^{'}\}$ be a useful ~system  in $\Phi$ . 
Then the elements of ${\cal T}_{\Delta}$ are called $\Delta-tableaux$,  the  
$J_{i}$ and $J_{i}^{'}$ are called  the $rows$ and the $columns$ of 
$\{J,J^{'}\}$ 
respectively.

\noindent
{\sc 3.4. Definition.}  Two $ \Delta$-tableaux $\bar{J}$ and $\bar{K}$  
are $row-equivalent$ ,  written $\bar{J}~\sim~\bar{K}$ ,  if there exists  $w \in {\cal W}(J)$ such that $\bar{K}~=~w~\bar{J}$.
The equivalence class which contains the $\Delta$-tableau $\bar{J}$ is 
$\{\bar{J}\}$ and is called a $\Delta-tabloid$.

Let $\tau_{\Delta}$ be the set of all $\Delta$-tabloids . 
It is clear that the number of distinct elements in 
$\tau_{\Delta}$ is  $[{\cal W}~:{\cal W}(J)]$  and by ( 2.7 )  we have
\[\tau_{\Delta}=\{ \{~d\bar{J}~\} \mid d\in D_{\Psi} \} \]
\noindent 
We note that if $\bar{J}=\{~J~;~J^{'}~\}$ then $dJ \subset \Phi^{+}$ but $dJ^{'}$ need not be a subset of $\Phi^{+}$ .

We now give an example to illustrate the construction of a $\Delta$-tabloid. 
In this example and later examples we use the following notation. 
If $\pi = \{\alpha_{1},\alpha_{2},...,\alpha_{n}\}$ is a simple system in $\Phi$ and $\alpha \in \Phi$, then $\alpha= {\displaystyle \sum_{i=1}^{n}a_{i}\alpha_{i}}$, where $a_{i} \in {\bf Z}$. From now on $\alpha$ is 
denoted by $a_{1}a_{2}...a_{n}$ and $\tau_{\alpha_{1}},
\tau_{\alpha_{2}},...,\tau_{\alpha_{n}}$ are denoted by 
$\tau_{1},\tau_{2},...,\tau_{n}$ respectively. 

\noindent
{\sc 3.5. Example.}  Let $\Phi={\bf D_{4}}$ 
 with simple system 
 \[\pi =\{\alpha_{1}=\epsilon_{1}-\epsilon_{2},
\alpha_{2}=\epsilon_{2}-\epsilon_{3},\alpha_{3}=\epsilon_{3}-\epsilon_{4},
\alpha_{4}=\epsilon_{3}+\epsilon_{4}\}\] 
\noindent
Let $\Psi_{1}={\bf A_{3}}$ be the subsystem 
of ${\bf D_{4}}$ with $J=\{1000,0100,0010\}$. Let $\Psi^{'}={\bf 2A_{1}}$ 
 be the subsystem of $\Phi$ which is 
contained in $\Phi \setminus \Psi$, with simple system $J^{'}=\{1101,0111\}$. 
Since ${\cal W}(J) \cap {\cal W}(J^{'})=<e>$ and   
${\cal W}(J^{\perp}) \cap {\cal W}(J^{'^{\perp}})=<e>$, 
then $\{J,J^{'}\}$ is a useful system in $\Phi$. 
Then $\tau_{\Delta}$ contains the $\Delta-tabloids$ ;

\begin {tabular}{lll}
$\{ \bar{J} \}$&=&$\{1000,0100,0010;1101,0111\}~$\\
$\{ \tau_{4}\bar{J} \}$&=&$\{1000,0101,0010;1100,0110 \}~$\\
$\{ \tau_{2}\tau_{4}\bar{J}\}$&=&$\{1100,0001,0110;1000,0010\}~$\\
$\{ \tau_{1}\tau_{2}\tau_{4}\bar{J} \}$&=&$\{0100,0001,1110;-1000,0010\}~$\\
$\{ \tau_{3}\tau_{2}\tau_{4}\bar{J} \}$&=&$\{1110,0001,0100;1000,-0010\}~$\\
$\{ \tau_{1}\tau_{3}\tau_{2}\tau_{4}\bar{J} \}$&=&$\{0110,0001,1100;-1000,-0010\}~$\\
$\{ \tau_{2}\tau_{1}\tau_{3}\tau_{2}\tau_{4}\bar{J} \}$&=&$\{0010,0101,1000;-1100,-0110\}~$\\
$\{ \tau_{4}\tau_{2} \tau_{1}\tau_{3}\tau_{2}\tau_{4}\bar{J} \}$&=&$\{0010,0100,1000;-1101,-0111\}$
\end{tabular}

\vspace{0.1cm}
The group ${\cal W}$ acts on $\tau_{\Delta}$ according to 
\[\sigma~\{\overline{wJ}\} =\{\overline{\sigma w J}\}~~~~~~for~all~\sigma \in {\cal W} .\]
\noindent
This action is well defined , for if $\{\overline{w_{1}J}\}=\{\overline{w_{2}J}\}$ , then there exists $\rho \in {\cal W}(w_{1}J)$ such that $\overline{\rho w_{1} J} =\overline{ w_{2}J}$ . 
Hence since $\sigma \rho \sigma^{-1} \in {\cal W}(\sigma w_{1}J)$ and $\overline{\sigma w_{2} J }=\overline{\sigma \rho w_{1} J} = (\sigma \rho \sigma^{-1} )(\overline{\sigma w_{1} J })$ , we have $\{\overline{\sigma w_{1}J}\}=\{\overline{\sigma w_{2}J}\}$ . 

Now if $K$ is arbitrary field , let $M^{\Delta}$ be the 
$K$-space whose basis elements are the $\Delta$-tabloids. Extend the action of ${\cal W}$ on $\tau_{\Delta}$ linearly on $M^{\Delta}$ , then $M^{\Delta}$ becomes  a $K{\cal W}$-module. Then we have the following lemma.

\noindent
{\sc 3.6. Lemma.} {\it The $K{\cal W}$-module $M^{\Delta}$  is a cyclic 
$K{\cal W}$-module generated by any one tabloid and 
$dim_{K}M^{\Delta} = [{\cal W}:{\cal W}(J)]$}

Now we proceed to consider the possibility of constructing a 
$K{\cal W}$-module which corresponds to the Specht module in the case of symmetric groups. In order to do this  we need to define a $\Delta$-polytabloid .

\noindent
{\sc 3.7. Definition.} Let $\{J,J^{'}\}$ be a useful system in $\Phi$ . Let
\begin {eqnarray*}
\kappa_{J^{'}}~=~\sum_{\sigma~\in~{\cal W}(J^{'})}~s~(~\sigma~)~\sigma~~~~~~and~~~e_{J,J^{'}}~=~\kappa_{J^{'}}~\{~\bar{J}~\}
\end{eqnarray*}
\noindent
where $s$ is the sign function defined in Section 2 . 
Then  $e_{J,J^{'}}$ is called the generalized $\Delta-polytabloid 
$ associated with $J$.

If $w \in {\cal W}(\Phi)$ , then 
\begin {eqnarray*}
w~\kappa_{J^{'}}&=&\sum_{\sigma~\in~{\cal W}(J^{'})}~s~(~\sigma~)~w~\sigma\\
&=&\sum_{\sigma~\in~{\cal W}(J^{'})}~s~(~\sigma~)~(w~\sigma~w^{-1})~w\\
&=&~\{~\sum_{\sigma~\in~{\cal W}(wJ^{'})}~s~(~\sigma~)~\sigma~\}~w
\end{eqnarray*}
\noindent
Hence, for all $w \in {\cal W}(\Phi)$, we have

\begin {tabular}{llll}
$w~e_{J,J^{'}}$&=&$\kappa_{wJ^{'}}~\{~\overline{wJ}~\}~=~e_{wJ,wJ^{'}}$~~~~~~~~~~ ~~(~3.1~)&
\end{tabular}

 Let $S^{J,J^{'}}$ be the subspace of $M^{\Delta}$ generated by $e_{wJ,wJ^{'}}$ where $w \in {\cal W}$. Then by (3.1) $S^{J,J^{'}}$ is a 
$K{\cal W}$-submodule of $M^{\Delta}$, 
which is called a $generalized$ $Specht$ $module$. Then we have the following theorem .

\noindent
{\sc 3.8. Theorem.} {\it The $K{\cal W}$-module $S^{J,J^{'}}$ is a cyclic submodule generated by any  $\Delta$-polytabloid.}

The following proposition notes some isomorphisms between Specht modules.

\noindent
{\sc 3.9. Proposition.} {\it Let $\{J,J^{'}\}$ be a useful system in $\Phi$ . Then we have the following isomorphisms:

\begin {tabular}{llll}
$(i)~If~w\in {\cal W}~,$&then&$S^{J,J^{'}}\cong S^{wJ,wJ^{'}}$&\\
$(ii)~If~w\in {\cal W}(J)~,$&then&$S^{J,J^{'}}\cong S^{J,wJ^{'}}$&\\
$(iii)~If~w\in {\cal W}(J^{'})~,$&then&$S^{J,J^{'}}\cong S^{wJ,J^{'}}$&
\end{tabular}}

Proposition 3.9 says  that a generalized Specht module is dependent only 
on the Dynkin diagram $\Delta$ and $\Delta^{'}$ of $J$ and $J^{'}$ 
respectively, thus, from now on it will be denoted by $S^{\Delta,\Delta^{'}}$.

A Specht module is spanned by the $e_{wJ,wJ^{'}}$ for all $ w \in {\cal W} $ ; the next lemma shows that we need only consider a certain subset of ${\cal W}$ .

\noindent
{\sc 3.10. Lemma.} {\it Let $\{J,J^{'}\}$ be a useful system in $\Phi$ . 
Then $S^{\Delta,\Delta^{'}}$ is generated by $e_{dJ,dJ^{'}}$ ,  
where $d \in D_{\Psi^{'}}$.}

\noindent
{\it Proof.} See Lemma 3.10 {\bf [8]}.

\noindent
{\sc 3.11. Lemma.} {\it Let $\{J,J^{'}\}$ be a useful system in $\Phi$ and 
let $d \in D_{\Psi}$. If $\{\overline{dJ} \}$ appears in $e_{J,J^{'}}$ then it appears only once.}

\noindent
{\it Proof.} See Lemma 3.11 {\bf [8]}.

\noindent
{\sc 3.12. Corollary.} {\it  If $\{J,J^{'}\}$ be a useful system in $\Phi$ , then $e_{J,J^{'}} \neq 0$ .}

The following lemma shows that the extra condition 
${\cal W}(J) \cap {\cal W}(J^{'})=<e>$  in our definition of a 
useful system is necessary. 
Unfortunately this condition which is a group theoretical one is 
not easily checked and it would be useful if it could be 
replaced by a criterion in terms of the root system only.

\noindent
{\sc 3.13. Lemma.} {\it If there exists $w \in {\cal W}(J) \cap 
{\cal W}(J^{'})$ such that $w$ has order 2 , and s($w$) = -1 then $e_{J,J^{'}}~=~0$.}

\noindent
{\it Proof.} See Lemma 3.13 {\bf [8]}.

\noindent
{\sc 3.14. Example.} Let $\Phi={\bf B_{3}}$ and 
$\pi=\{\alpha_{1}=\epsilon_{1}-\epsilon_{2},\alpha_{2}
=\epsilon_{2}-\epsilon_{3},\alpha_{3}=\epsilon_{3} \}$. Let $\Psi=
{\bf 3A_{1}}$ be the subsystem of $\Phi$ with simple system 
$J=\{\alpha_{1}=\epsilon_{1}-\epsilon_{2},\tilde{\alpha}=\epsilon_{1}+\epsilon_{2},\alpha_{3}=\epsilon_{3} \}$ and 
let $\Psi^{'}={\bf 3A_{1}}$ be the subsystem of $\Phi$ with $J^{'}=\{\alpha_{2}=\epsilon_{2}-\epsilon_{3},\alpha_{1}+\alpha_{2}+\alpha_{3}=\epsilon_{1},\alpha_{2}+2~\alpha_{3}=\epsilon_{2}+\epsilon_{3}\}$. Then $\Psi \cap \Psi^{'}=\emptyset$. But
\newline
${\cal W}(J)~~$=$~\{e,\tau_{1},\tau_{3},\tau_{1}\tau_{3},
\tau_{2}\tau_{3}\tau_{1}\tau_{2}\tau_{3}\tau_{1}\tau_{2},\tau_{2}\tau_{3}\tau_{1}\tau_{2}\tau_{3}\tau_{1}\tau_{2}\tau_{1},
\tau_{3}\tau_{2}\tau_{3}\tau_{1}\tau_{2}\tau_{3}\tau_{1}\tau_{2},\tau_{3}\tau_{2}\tau_{3}\tau_{1}\tau_{2}\tau_{3}\tau_{1}\tau_{2}\tau_{1}\}$ 
\newline
${\cal W}(J^{'})~$=$~\{e,\tau_{2},
\tau_{1}\tau_{2}\tau_{3}\tau_{2}\tau_{1},\tau_{3}\tau_{2}\tau_{3},\tau_{3}\tau_{2}\tau_{3}\tau_{1}\tau_{2}\tau_{3}\tau_{1}\tau_{2}\tau_{1},\tau_{3}\tau_{2}\tau_{3}\tau_{2},
\tau_{3}\tau_{2}\tau_{3}\tau_{1}\tau_{2}\tau_{3}\tau_{2}\tau_{1},\tau_{1}\tau_{2}\tau_{3}\tau_{1}\tau_{2}\tau_{1}\}$ 

\noindent 
It follows that $w=~\tau_{3}\tau_{2}\tau_{3}\tau_{1}\tau_{2}\tau_{3}\tau_{1}\tau_{2}
\tau_{1}~ \in {\cal W}(J) \cap {\cal W}(J^{'}) $ and $e_{J,J^{'}}=0$.

\noindent
{\sc 3.15. Lemma.} {\it  Let $\{J,J_{1}^{'}\}$ and $\{J,J_{2}^{'}\}$ be useful 
systems in $\Phi$ . If $\Psi_{1}^{'} \subseteq \Psi_{2}^{'}$ , 
then $S^{J,J_{2}^{'}}$ is a $K{\cal W}$-submodule of $S^{J,J_{1}^{'}}$ , 
where $J_{1}^{'}$ and $J_{2}^{'}$  are simple systems for $\Psi_{1}^{'}$ 
and $\Psi_{2}^{'}$ respectively.}

Now we  consider under what conditions  
$S^{\Delta,\Delta^{'}}$ is irreducible .

\noindent
{\sc 3.16. Lemma.} {\it Let $\{J,J^{'}\}$ be a useful system in $\Phi$ and  let $ d \in D_{\Psi}$ . Then the following conditions are equivalent:
\newline
(i) $\{~\overline{dJ}~\}$ appears with non-zero coefficient in $e_{J,J^{'}}$
\newline
(ii) There exists $\sigma \in {\cal W}(J^{'})$ such that $\sigma \{~\bar{J}~\}~=~\{\overline{dJ}~\}$
\newline
(iii) There exists $\rho \in {\cal W}(J)$  and $\sigma \in {\cal W}(J^{'})$ such that $d~=~\sigma~\rho$}

\noindent
{\it Proof.}  See Lemma 3.16 {\bf [8]}.

\noindent
{\sc 3.17. Lemma.} {\it Let $\{J,J^{'}\}$ be a useful system in $\Phi$ and  let $ d \in D_{\Psi}$ . If $\{~\overline{dJ}~\}$ appears in
 $e_{J,J^{'}}$ then $d~\Psi \cap \Psi^{'}~=~\emptyset $.}

\noindent
{\it Proof.} See Lemma 3.17 {\bf [8]}.

\noindent
{\sc 3.18. Lemma.} {\it Let $\{J,J^{'}\}$  be  a useful system in $\Phi$ and 
let $ d \in D_{\Psi}$. Let $d~\Psi \cap \Psi^{'} \not = \emptyset$. 
Then $\kappa_{J^{'}}\{~\overline{dJ}~\} = 0 $.}

The converse of Lemma 3.17 is not true in general, 
which leads to the following definition .

\noindent
{\sc 3.19. Definition.}  A useful system $\{J,J^{'}\}$ in $\Phi$ is 
called a $good~system$  if $d~\Psi \cap \Psi^{'} = \emptyset$ for $d \in D_{\Psi}$ then $\{~\overline{dJ}~\}$ 
appears with non-zero coefficient in $e_{J,J^{'}}$.

\noindent
{\sc 3.20. Lemma.} {\it  Let $\{J,J^{'}\}$ be a good system in $\Phi$ and  let $d \in D_{\Psi}$ .
\newline
(i) If $\{\overline{dJ}\}$ does not appear in $e_{J,J^{'}}$  then $\kappa_{J^{'}}\{~\overline{dJ}~\} = 0 $.
\newline
(ii) If $\{\overline{dJ}\}$  appears in $e_{J,J^{'}}$  then there exists $ \sigma \in {\cal W}(J^{'})$ such that}
\[\kappa_{J^{'}}\{~\overline{dJ}~\}~=~s~(~\sigma~)~e_{J,J^{'}} \]

\noindent
{\it Proof.} See Lemma 3.20 {\bf [8]}.

\noindent
{\sc 3.21. Corollary.} {\it  Let $\{J,J^{'}\}$ be a 
good system in $\Phi$.  If $m \in M^{\Delta}$ then $\kappa_{J^{'}} 
m$ is a
 multiple of $e_{J,J^{'}}$.}

We now define a bilinear form $ <~,~>$ on 
$M^{\Delta}$ by setting

\[
<~\{\bar{J_{1}}\}~,~\{\bar{J_{2}}\}~>~=
\left\{
\begin{array}{ll}
1& \mbox{if~$\{\bar{J_{1}}\}~=~\{\bar{J_{2}}\}$}\\
0& \mbox{otherwise} 
\end{array}
\right.
\]

\noindent
This is a symmetric, non-singular, ${\cal W}$-invariant bilinear form
 on $M^{\Delta}$ .

Now we can prove James' submodule theorem in this general setting.

\noindent
{\sc 3.22. Theorem.} {\it  Let $\{J,J^{'}\}$ be a good 
system in $\Phi$.  Let U be submodule of $M^{\Delta}$ .
 Then either $S^{\Delta , \Delta^{'}} \subseteq U$ or $U 
\subseteq S^{\Delta ,
 \Delta^{'^\perp}}$ , where $S^{\Delta,\Delta^{'^\perp}}$
 is the complement of $S^{\Delta,\Delta^{'}}$ in  $M^{\Delta} .$}

\noindent
{\it Proof.} See Theorem 3.22 {\bf [8]}.

We can now prove our principal result.

\noindent
{\sc 3.23. Theorem.} {\it Let $\{J,J^{'}\}$ be a
 good system in $\Phi$ . The $K{\cal W}$-module 
\newline
$D^{\Delta,\Delta^{'}}~=~S^{\Delta,\Delta^{'}}~
/~S^{\Delta,\Delta^{'}}\cap S^{\Delta,\Delta^{'^\perp}}$ 
is zero or irreducible.}

\noindent
{\it Proof.} 
 If $U$ is a submodule of $S^{\Delta,\Delta^{'}}$ 
then $U$ is a 
submodule of $M^{\Delta}$ and by Theorem 3.22 
either $S^{\Delta,\Delta^{'}} 
\subseteq U$ in which case $U~=~
S^{\Delta,\Delta^{'}} $ or $U \subseteq 
S^{\Delta,\Delta^{'^\perp}}$ and $U
 \subseteq ~S^{\Delta,\Delta^{'}}~\cap ~
S^{\Delta,\Delta^{'^\perp}}$ , which 
completes the proof.

In the case of K = {\bf Q} or any field of
 characteristic zero $<~,~>$ is an inner product and
 $D^{\Delta,\Delta^{'}}~=~S^{\Delta,\Delta^{'}}$.
 Thus if for a 
subsystem $\Psi$ of $\Phi$
 a good system $\{ J,J^{'}\}$ can be found , 
then we have a construction for irreducible 
$K{\cal W}$-modules . Hence it is essential to 
show for each subsystem 
that a good system exists which satisfies
 Definition 3.19.

In the following example, we show how a good system 
may be constructed in all cases for the finite reflection group of 
type ${\bf G_{2}}$. In {\bf [9]}, we present an algorithm 
for constructing a good system for certain subsystems; 
indeed this
 algorithm will give a good system with additional
 properties which will lead to the construction
 of a $K$-basis for our Specht
 modules  $S^{\Delta,\Delta^{'}}$, which correspond to the basis consisting
 of standard tableaux in the case of symmetric groups.

\noindent
{\sc 3.24. Example.} Let $\Phi={\bf G_{2}}$ with simple system
$\pi = \{\alpha_{1}=\epsilon_{1}-\epsilon_{2},\alpha_{2}
=-2\epsilon_{1}+\epsilon_{2}+\epsilon_{3}\}$. 
Let $g_{1}=e,g_{2}=\tau_{2},g_{3}=\tau_{1}\tau_{2},
g_{4}=(\tau_{1}\tau_{2})^{2},
g_{5}=(\tau_{1}\tau_{2})^{3},g_{6}=\tau_{1}$ 
be representatives of conjugate classses  
$C_{1},C_{2},C_{3},C_{4},C_{5},C_{6}$ respectively of ${\cal W}({\bf G_{2}})$. 
Then the character table of ${\cal W}({\bf G_{2}}$) is

\vspace{0.5cm}
\begin{tabular}{r|rrrrrr}
&$C_{1}$&$C_{2}$&$C_{3}$&$C_{4}$&$C_{5}$&$C_{6}$\\
\hline  
$\chi_{1}$&1&1&1&1&1&1\\ 
$\chi_{2}$&1&-1&1&1&1&-1\\
$\chi_{3}$&1&-1&-1&1&-1&1\\
$\chi_{4}$&1&1&-1&1&-1&-1\\
$\chi_{5}$&2&0&-1&-1&2&0\\
$\chi_{6}$&2&0&1&-1&-2&0
\end{tabular}  

\newpage
The non-conjugate  subsystems of ${\bf G_{2}}$ are:

\begin{tabular}{lllll}
(1)  $\Psi_{1}$&=&${\bf A_{2}}$& 
 with simple system~ $J_{1}~= $&$~\{01,31\}~$
\\
(2)  $\Psi_{2}$&=&${\bf A_{1}+\tilde{A_{1}}}$&  
with simple system ~$J_{2}~=$&$~\{10,32\}~$
\\
(3)  $\Psi_{3}$&=&${\bf A_{1}}$ & 
with simple system~ $J_{3}~=$&$~\{10\}~$
\\
(4)  $\Psi_{4}$&=&${\bf \tilde{A_{1}}}$& 
 with simple system ~$J_{4}~=$&$~\{01\}$
\\
(5)  $\Psi_{5}$&=&$\emptyset$ &
 with simple system ~$J_{5}~=$&$~\emptyset$
\\
(6)  $\Psi_{6}$&=&${\bf G_{2}}$ &
 with simple system~ $J_{6}~=$&$~\{10,01\}$
\end{tabular}
\vspace{0.5cm}

Let $\Psi_{4}={\bf \tilde{A_{1}}}$ be the subsystem of $\Phi$ 
with simple system $J_{4}=\{01\}$. 
Let $\Psi_{1}^{'}={\bf A_{1}}+{\bf \tilde{A_{1}}}$ be the subsystem of $\Phi$
 which is contained in $\Phi \setminus \Psi_{4}$, with simple system 
 $J_{1}^{'}=\{11,31\}$. Since ${\cal W}(J_{4}) \cap 
 {\cal W}(J_{1}^{'})=<e>$ and   
${\cal W}(J_{4}^{\perp}) \cap {\cal W}(J_{1}^{'^{\perp}})=<e>$, 
then $\{J_{4},J_{1}^{'}\}$ is a useful system in $\Phi$. 
Then $\tau_{\Delta_{4}}$ 
contains $\Delta_{4}$-tabloids:
\[\{\overline{J_{4}}\}=\{01;11,31\},\{\overline{\tau_{1}J_{4}}\}=\{31;21,01\} \]
\[\{\overline{\tau_{2}\tau_{1}J_{4}}\}=\{32;21,-01\},
\{\overline{\tau_{1}\tau_{2}\tau_{1}J_{4}}\}=\{32;11,-31\}\]
\[\{\overline{\tau_{2}\tau_{1}\tau_{2}\tau_{1}J_{4}}\}=\{31;10,-32\},
 \{\overline{\tau_{1}\tau_{2}\tau_{1}\tau_{2}\tau_{1}J_{4}}\}=\{01;-10,-32\}\]

For $d=e, \tau_{2}\tau_{1},\tau_{1}\tau_{2}\tau_{1}, 
\tau_{1}\tau_{2}\tau_{1}\tau_{2}\tau_{1}$
  we have $d\Psi_{4} \cap \Psi_{1}^{'} = \emptyset$. Since
  
\[e_{J_{4},J_{1}^{'}}=\{\bar{J_{4}}\}-\{\overline{
\tau_{2}\tau_{1}J_{4}}\}-
  \{\overline{\tau_{1}\tau_{2}\tau_{1}J_{4}}\}+
  \{\overline{\tau_{1}\tau_{2}\tau_{1}\tau_{2}\tau_{1}J_{4}}\} \]
 \noindent
 then $\{J_{4},J_{1}^{'}\}$ is  a good system in $\Phi$. 

Now let $K$ be a field and Char$K$ = 0 . Let 
$M^{\Delta_{4}}$ be $K$-space whose basis are the $\Delta_{4}$-tabloids. 
Let $S^{\Delta_{4},\Delta_{1}^{'}}$ be the corresponding 
 $K{\cal W}$-submodule of 
 $M^{\Delta_{4}}$, then  by definition of 
 the Specht module we have 
\[S^{\Delta_{4},\Delta_{1}^{'}}=
Sp~\{~e_{J_{4},J_{1}^{'}}~,~e_{\tau_{1}J_{4},\tau_{1}J_{1}^{'}}~\}\]

\noindent
where

$\begin{array}{lll}
e_{J_{4},J_{1}^{'}}&=&
\{\bar{J_{4}} \}
-\{\overline{ \tau_{2}\tau_{1}J_{4}} \}
-\{\overline{ \tau_{1}\tau_{2}\tau_{1}J_{4}} \}
+\{\overline{ \tau_{1}\tau_{2}\tau_{1}\tau_{2}\tau_{1}J_{4}} \}
\\
e_{\tau_{1}J_{4},\tau_{1}J_{1}^{'}}&
=&
\{\overline{\tau_{1}J_{4}} \}
-\{ \overline{\tau_{2}\tau_{1}J_{4}} \}
-\{\overline{ \tau_{1}\tau_{2}\tau_{1}J_{4}} \}
+\{\overline{ \tau_{2}\tau_{1}\tau_{2}\tau_{1}J_{4}} \}
\end{array}$

Let $T_{4}^{1}$ be the matrix representation of ${\cal W}$ afforded by 
$S^{\Delta_{4},\Delta_{1}^{'}}$ with character $\psi_{4}^{1}$  and 
let $\tau_{1}\tau_{2}$ be the representative of the
conjugate class $C_{3}$. Then

{\setlength{\arraycolsep}{0.1cm}
$\begin{array}{ll}
\tau_{1}\tau_{2}(e_{J_{4},J_{1}^{'}})&=
e_{\tau_{1}J_{4},\tau_{1}J_{1}^{'}}~-~e_{J_{4},J_{1}^{'}}\\
\tau_{1}\tau_{2}(e_{\tau_{1}J_{4},\tau_{1}J_{1}^{'}})&
=~-~e_{J_{4},J_{1}^{'}}
\end{array}$}

Thus we have

$T_{4}^{1}~(~\tau_{1}\tau_{2}~)~=
~\left( \begin{array}{cc}
-1&1\\
-1&0
\end{array}
\right) $ and $\psi_{4}^{1}(\tau_{1}\tau_{2})=-1$.

By a similar calculation to the above it can be shown that 
$\psi_{4}^{1}=\chi_{5}$. By the same method to the above, we have

\vspace{0.2cm}

{\setlength{\arraycolsep}{0.2cm}
$\begin{array}{l|l|l|l}
\Psi&\Psi^{'}&~~~J^{'}&Ch\\
\hline
{\bf A_{2}}&{\bf
A_{1}}&\{10\}&\chi_{4}
\\
{\bf A_{1}}+{\bf \tilde{A_{1}}}&
{\bf \tilde{A_{1}}}&\{01\}&\chi_{5}
\\
{\bf A_{1}}&{\bf
A_{2}}&\{01,31\}&\chi_{3}
\\
{\bf \tilde{A_{1}}}&
{\bf A_{1}}+{\bf \tilde{A_{1}}}&\{11,31\}&\chi_{5}
\\
{\bf G_{2}}&
~~\emptyset&~~\emptyset&\chi_{1}
\\
~~\emptyset&{\bf G_{2}}&\{10,01\}&\chi_{2}

\end{array}$}
\vspace{0.2cm}

We note that the irreducible modules corresponding  to the characters
$\chi_{6}$ have not been obtained. We now show how
an additional irreducible character is obtained. Let $\Psi_{2}^{'}={\bf
A_{1}}$ be the subsystem of $\Phi$ with simple system
$J_{2}^{'}=\{11\}$. Then $\{J_{4},J_{2}^{'}\}$ is a useful system in $\Phi$.
Since $\Psi_{2}^{'} \subset  \Psi_{1}^{'}$, by Lemma 3.15 
$S^{\Delta_{1},\Delta_{1}^{'}}$ is a $K{\cal W}$-submodule of 
$S^{\Delta_{1},\Delta_{2}^{'}}$. By a similar calculation, 
the corresponding character of ${\cal W}$ afforded by 
$S^{\Delta_{1},\Delta_{2}^{'}}/S^{\Delta_{1},\Delta_{1}^{'}}$ is $\chi_{6}$.
Thus we have obtained a complete set of irreducible modules for ${\bf G_{2}}$.

\begin{center}        
REFERENCES
\end{center}

{\bf 1.}  E. Al-Aamily , A . O .  Morris and  M . H . Peel,
The~ representations~ of~ the~ Weyl groups~ of ~type ~$B_{n}$, 
$Journal~ of~ Algebra$ {\bf 68} ( 1981), 298--305.

{\bf 2.}  A . Borel and  J de Siebenthal,
Les sous-groupes  ferm$\grave{e}$s connexes de rang maximium des 
groupes de Lie clos,  $Comment.~  Math.~ Helv.$, {\bf 23}(1949), 200--221.

{\bf 3.}  N . Bourbaki,
{\em Groupes et alg$\grave{e}$bres de Lie, Chapters 4,5,6},  
Actualites Sci. Induct  1337, (Hermann, Paris, 1968)

{\bf 4.}  R. W. Carter, 
Conjugacy~ classes~ in ~the Weyl~ group,
$Comp.~  Math.$, {\bf 25} (1972), 1--59.

{\bf 5.}  R . W . Carter, 
{\em Simple Groups of Lie Type},  
(Wiley, London,  Newyork, Sydney, Toronto, 1989)

{\bf 6.}  E . B . Dynkin, 
Semisimple subalgebras of semisimple Lie algebras,
 {\it Amer. Math.  Soc. Trans.(2)}  {\bf 6}, (1957),  111--244.

{\bf 7.}  L . C . Grove and C . T . Benson ,
{\em Finite Reflection Groups}, (Springer Verlag, Newyork , Berlin , 
Heilderberg , Tokyo , 1985 ) . 

{\bf 8.}  S. Hal\i c\i o{\u g}lu and  A. O. Morris, 
Specht~ Modules~ for~ Weyl~ Groups, \\
$Contributions$ $to~ Algebra$ $and~ Geometry$, {\bf 34} (1993), 257--276.

{\bf 9.}  S. Hal\i c\i o{\u g}lu, 
 A Basis of Specht~ Modules~ for~ Weyl~ Groups, 
submitted~ for~ publication.

{\bf 10.} J . E . Humphreys,
{\em  Introduction to Lie algebras and representation theory},  
Graduate Texts in Mathematics,  Volume 9 ( Springer-Verlag, Berlin, 1972 )

{\bf 11.} J . E . Humphreys,
{\em Reflection Groups and Coxeter Groups} ,  
 ( Cambridge University Press , Cambridge , 1990)
         
{\bf 12.} G. D. James,
The~ irreducible~ representations~ of~ the~ symmetric~ group, \\
$Bull.~ Lond.~ Math.~ Soc$, {\bf 8} (1976), 229--232.

{\bf 13.}  G . D . James , A . Kerber,
{\em The Representation Theory of the Symmetric Group},  
Addison-Wesley Publishing Company ( London, 1981 ).

{\bf 14.} A. O. Morris,
Representations~ of~ Weyl~ groups~ over~ an~ arbitrary~ field, 
$Ast\grave{e}risque$ {\bf 87-88} (1981), 267--287.

{\bf 15.} A. O. Morris and A. J. Idowu, 
Some combinatorial results for Weyl groups,  
{\it Proc. Camb. Phil. Soc.}  {\bf 101} (1987), 405--420

\noindent
{\sc Department of Mathematics\\ 
Ankara University\\ 
06100~  Tando{\u g}an\\
Ankara Turkey}\\

\end{document}